\documentclass{article}

\usepackage{amssymb,dsfont,latexsym}

\newtheorem{thm}{Theorem}

\newtheorem{prop}[thm]{Proposition}
\newtheorem{lem}[thm]{Lemma}
\newtheorem{cor}[thm]{Corollary}

\newcommand\enu[1]{\smallskip\newline\makebox[5mm][l]{\rm(#1)}}

\newcommand\bp{\noindent{\it Proof.}\ }

\newcommand\2{{\frac{1}{2}}}

\begin{document}

\author{Erling St{\o}rmer}

\date{28-3-2008 }

\title{Separable states and positive maps II}

\maketitle
\begin{abstract}
Using the natural duality between linear functionals on tensor
products of C*-algebras with the trace class operators on a
Hilbert space $H$ and linear maps of the C*-algebra into $B(H)$,
we give two characterizations of separability, one relating it to
abelianness of the definite set of the map, and one on tensor
products of nuclear and UHF $C^*-$algebras.

\end{abstract}

\section*{Introduction}

The present paper is a continuation of \cite{st5}. We continue to
study the duality of positive linear maps from a $C^*-$algebra $A$
into the bounded operators on a Hilbert space and linear
functionals on the tensor product of $A$ with the trace class
operators. The paper consists of two rather independent parts.  In
the first we consider the definite set $D_\phi$ of a map $\phi$,
i.e. the Jordan subalgebra of $A$ on which $\phi$ restricts to a
Jordan homomorphism. The main result states that if
$\phi(A)=\phi(D_\phi)$, then the dual functional of $\phi$ is a
separable positive functional if and only if $\phi(A)$ is an
abelian $C^*-$algebra. This is applied to give necessary and
sufficient conditions for the dual functional of the canonical
projection of a von Neumann algebra onto the fixed point set of a
positive unital map to be separable.

In the second part we prove a variation of the Horodecki theorem
\cite{3 Hor} characterizing separable states by positive maps in
the setting of tensor products of a nuclear $C^*-$algebra $A$ and
an UHF-algebra $B$. We show that a state $\rho$ on $A\otimes B$ is
the $w^*$-limit of convex sums of product states if and only if
the compositions $\rho \circ (\iota\otimes\psi)$ are positive for
all positive maps $\psi$ of $B$ into itself, $\iota$ being the
identity map on $A$.

For the reader's convenience we recall some of the concepts
considered in the paper.

Let $A$ be a $C^*-$algebra, or just an operator system. Let $\1T$
denote the trace class operators on a Hilbert space $H$,
identified with $B(H)$ when $H$ is finite dimensional. If $B$ is
another $C^*-$algebra let $B(A,B)$ (resp.$B(A,B)^+$) denote the
bounded linear (resp. positive linear) maps of $A$ into $B$. We
also use the notation $B(A,H)$ for $B(A,B(H))$. For the following
see \cite {st5}. If $\phi\in B(A,B)$ its dual functional
$\tilde{\phi}$ on the projective tensor product
$A\widehat{\otimes}\1T$ is defined by
$$
\tilde{\phi}(a\otimes b)=Tr(\phi(a)b^t);
$$
$Tr$ is the usual trace on $B(H)$, and $b^t$ is the transpose of
$b$ with respect to a given orthonormal basis for $H$. A
\textit{mapping cone} is a cone $0\neq K \subseteq B(B(H),H)^+$
which is closed in the topology of bounded pointwise weak
convergence, and which has a dense set of maps $\alpha$ such that
for all $a,b\in B(H)$ the map $x\mapsto a\alpha(bxb^*)a^*$ belongs
to $K$. We say $\phi\in B(A,H)^+$ is \textit{K-positive} if
$\tilde\phi$ is positive on the cone
$$
P(A,K)=\{x\in A\widehat{\otimes}\1T: \iota\otimes\alpha(x)\geq 0
  \ \forall\ \alpha\in K\}.
$$
Note that if $x\in P(A,K)$ and $\beta\in B(A,A)$ is completely
positive, then
$$
(\iota\otimes\alpha)(\beta\otimes\iota(x))=(\beta\otimes\iota)\circ(\iota\otimes\alpha)(x)\geq
0,
$$
hence $\beta\otimes\iota(x)\in P(A,K)$. It follows that $x\in
P(A,K)$ if and only if $x\in A\widehat{\otimes}\1T$ and
$\beta\otimes\alpha(x)\geq 0 $ for all $\alpha\in K, \beta\in
B(A,A)$ completely positive.

Following \cite{HSR} we say a map $\phi\in B(A,H)^+$ is
\textit{entanglement breaking} if it is $S(H)-$positive, where
$S(H)$ is the mapping cone generated by maps of the form
$$
\phi(x)=\sum_{i=1}^{n}\omega_{i}(x)a_i,
$$
where $\omega_i$ is a normal state on $B(H)$ and $a_i \in B(H)^+$.
By \cite{st5},Theorem 2,see also \cite{HSR}, $\phi$ is
entanglement breaking if and only if $\tilde{\phi}$ is a separable
positive linear functional.

A linear set $J$ of self-adjoint operators is called a
\textit{Jordan algebra} if $a,b\in J$ implies their Jordan product
$a\circ b= \2(ab+ba)\in J.$ $J$ is called a JC-algebra (resp.
JW-algebra) if it is norm (resp. weakly) closed. In particular the
self-adjoint part $A_{sa}$ of a $C^*-$algebra (resp. von Neumann
algebra) is a JC-algebra (resp. JW-algebra).

If $\phi\in B(A,B)^+$ with $A,B$ $C^*-$algebras, the
\textit{definite set $D_\phi$} of $\phi$ is the set \textmd{of}
operators $a\in A$ such that $\phi(a^*a)=\phi(a)^* \phi(a).$ In
particular the self-adjoint part $(D_\phi)_{sa}$ is the set of
self-adjoint operators such that $\phi(a^2)=\phi(a)^2.$ Then
$(D_\phi)_{sa}$ is a JC-algebra. If $a\in D_\phi, b\in A$ then
$\phi(a\circ b)=\phi(a)\circ\phi(b),$ see \cite{st4}. Note that in
\cite{st5} we only considered the self-adjoint part of $D_\phi$
and called that the definite set.

The author is indebted to M.B.Ruskai for some valuable comments.

\section{Separability and abelianness}

It was shown in \cite{st5} that a map in $B(A,H)^+$ of the form
$\phi(x)=\sum\omega_i(x)a_i$ with $\omega_i$ states of $A$ and
$a_i\in B(H)^+$ has definite set whose image by $\phi$ is an
abelian $C^*-$algebra. In this section we shall elaborate on this
and show that if $\phi\in B(A,H)^+$ has image
$\phi(A)=\phi(D_\phi)$ then $\tilde{\phi}$ is separable if and
only if $\phi(A)$ is abelian.  In the von Neumann algebra case
when $\phi\in B(M,M)^+ $ is normal we apply this to the positive
projection map $P_\phi$ onto the fixed point set $M_\phi$ of
$\phi$ in M. But first we show some preliminary results.  The
first is a simplified dual characterization of decomposable and
entanglement breaking maps.

\begin{prop}\label{prop1}
Let $A$ be an operator system, $H$ a Hilbert space and $\phi\in
B(A,H)^+$. Let $t$ denote the transpose map on $B(H).$ Then we
have
 \enu{i}$\phi$ is decomposable if and only if $\tilde{\phi}$ is
 positive on the cone
 $$
 \{x\in (A\widehat{\otimes}\1T)^+: (\iota\otimes
 t)(x)\geq 0\}.
 $$
 \enu{ii} $\phi$ is entanglement breaking if and only if
 $\tilde{\phi}$ is positive on the cone
 $$
 \{x\in A\widehat{\otimes}\1T: (\iota\otimes\omega)(x)\geq 0
  \ \forall\hbox{ normal states } \omega \hbox{ of } B(H)\}.
 $$
 \end{prop}
 \bp
Ad (i).Let $K$ be the mapping cone of all decomposable maps in
 $B(B(H),H).$ Then $\phi$ is $ K$-positive if and only if
 $\tilde{\phi}$ is positive on the cone
 $$
 \{x\in A\widehat{\otimes}\1T: (\iota\otimes\alpha)(x)\geq 0
  \ \forallÊ\ \alpha\in K\},
 $$
 or, since each copositive map is of the form $\beta\circ t$ with
 $\beta$ completely positive, it follows that if $ \iota\otimes t)(x)\geq
 0$ then $ (\iota\otimes \beta\circ t)(x)= (\iota\otimes \beta)\circ(\iota\otimes t)(x)\geq 0.$
 Thus $\phi$ is $K$-positive if and only if $\tilde{\phi}$ is
 positive  on the cone $\{x\in (A\widehat{\otimes}\1T)^+: (\iota\otimes
 t)(x)\geq 0\}$. But then it follows from \cite{st3} Theorem 3.6,
 that $\phi$ is $K$-positive if and only if $ \phi$ belongs to the
 closed cone generated by all maps $\alpha\circ\psi, \alpha\in K, \psi\in
 B(A,H)$ completely positive, in other words, if and only if
 $\phi$ is decomposable, completing the proof of (i).

 Ad (ii). Assume $\tilde{\phi}$ is positive on the cone
 $$
 \{x\in A\widehat{\otimes}\1T:(\iota\otimes\omega)(x)\geq 0
  \ \forall \ \hbox{normal states }\omega\ \hbox{of } B(H)\}.
 $$
 If $(\iota\otimes\omega)(x)\geq 0$, then
 $$
(\iota\otimes a\omega)(x)= (1\otimes a)(\iota\otimes\omega)(x)\geq
0
$$
if $a\geq 0,$ since $(1\otimes a)$ and $(\iota\otimes\omega)(x)$
commute. Thus $(\iota\otimes\alpha)(x)\geq 0$ for all $\alpha\in
S(H),$ so $\phi$ is $S(H)-$positive, hence entanglement breaking, 
by \cite{st5}Theorem 2. The proof is complete.

\begin{lem}\label{lem1}
Let $A$ be a unital $C^*-$algebra and $\phi\in B(A,H)$ unital and
completely positive. Then the definite set $D_\phi$ is a
$C^*-$subalgebra of $A$.
\end{lem}

\bp Let $\phi=V^*\pi V $ be the Stinespring dilation of $\phi,$
where $\pi$ is a representation of $A$ on a Hilbert space $K,$ and
$V$ is a linear map of $H$ into $K.$ Since $\phi$ is unital,
$V^*V=1,$ so $V$ is an isometry. Let $p=VV^*.$ Then $p$ is a
projection in $B(K)$ such that the map $p\pi p$ is the composition
of $\phi$ with the isomorphism $x\mapsto VxV^*$ of $B(H)$ into
$B(K).$ Thus $x\in D_\phi$ if and only if $x\in D_{p\pi p}.$ But
then $x\in (D_\phi)_{sa}$ if and only if
$$
p\pi(x)p\pi(x)p=(p\pi(x)p)^2=p\pi(x^2)p=p\pi(x)^2 p,
$$
hence if and only if $p\pi(x)(1-p)\pi(x)p=0,$ i.e.
$p\pi(x)(1-p)=0$, hence if and only if
$p\pi(x)=p\pi(x)p=(p\pi(x)p)^*=\pi(x)p.$ Thus $x\in D_\phi$ if and
only if $p\pi(x)=\pi(x)p.$ Thus $x,y\in D_\phi$ implies their
product $xy\in D_\phi,$ proving the lemma.

\medskip

The next lemma is closely related to Theorem 10 in \cite{st5}

\begin{lem}\label{lem2}
Let $\phi\in B(A,H)$ be unital and entanglement breaking. Then
$\phi(D_\phi)$ is an abelian $C^*-$algebra.
\end{lem}
\bp Since an entanglement breaking map is completely positive,
$D_\phi$ is a $C^*-$algebra by Lemma 2. Furthermore $\phi$ is a
Jordan homomorphism on $D_\phi.$ By \cite{st1} there are two
orthogonal central projections  $e$ and $f$ with sum 1 in the von
Neumann algebra generated by $\phi(D_\phi)$ such that
$\phi_1(a)=\phi(a)e $ is a *-homomorphism, and
$\phi_2(a)=\phi(a)f$ is a *-antihomomorphism. Since $\phi$ is
entanglement breaking the composition of $\phi$ with the transpose
map is also entanglement breaking, hence completely positive. Thus
both $\phi_1$ and $\phi_2$ are both homomorphisms and
antihomomorphisms. But this is possible only if their images are
abelian. The proof is complete.

\medskip
By the above together with some results in \cite{st5} we can now
conclude the following.

\begin{thm}\label{thm1}
Let $A$ be a unital $C^*-$algebra and $\phi$ a unital map in
$B(A,H)^+$ such that $\phi(A)=\phi(D_\phi).$ Then the following
conditions are equivalent:
 \enu{i} $\phi$ is entanglement breaking.
 \enu{ii} $\phi(A)$ is an abelian $C^*-$algebra.
 \enu{iii} $\tilde{\phi}$ is separable.
 \end{thm}
 \bp
 By Theorem 2 in \cite{st5} $(i)\Leftrightarrow(iii)$. By Corollary 3
 in \cite{st5} $(ii)\Rightarrow(iii)$, and by Lemma ~\ref{lem2}
 $(i)\Rightarrow(ii)$. The proof is complete.

\medskip
 Let $M\subset B(H)$ be a von Neumann algebra, and let $\phi\in
 B(M,M)^+$. Then $\phi^n \in B(M,M)^+$ for all $n\in \0N$. It
 was shown in the proof of \cite{ES} Corollary 1.6, that if $\phi$
 is normal, then then the maps $\phi_n =\frac{1}{n}(\phi + \phi^2
 +...+\phi^n)$ converge in the point-ultraweak topology to a
 positive unital projection map $P_\phi$, so $P_\phi =(P_\phi)^2$,
 of $M$ onto the fixed point set $M_\phi =\{x\in M: \phi(x)=x\}$
 of $\phi.$ We call this projection map the \textit{averaging
 projection} for $\phi.$ Suppose furthermore that there exists a
 faithful normal state $\rho$ on $M$ such that $\rho\circ\phi =
 \rho.$ Then $\rho\circ\phi^n = \rho$ for all n, hence $\rho\circ P_\phi
 =\rho$. It follows that $P_\phi$ is faithful and normal. Indeed,
 if $(a_\alpha)$ is an increasing net in $M$, and $a_\alpha \nearrow
 a,$ then
 $$
 \rho(P_{\phi}(a))=\rho(a)=\lim \rho(a_\alpha)=\lim \rho
 P_{\phi}(a_\alpha)\leq\rho P_{\phi}(a),
 $$
 so by faithfulness of $\rho, P_{\phi}(a) = \lim
 P_{\phi}(a_\alpha)$, so $P_\phi$ is normal. Similarly one shows
 that $P_\phi$ is faithful. Then by \cite{ES} Corollary 1.5,
 $(M_\phi)_{sa}=P_{\phi}(M_{sa})$ is a JW-algebra. Furthermore $M_\phi
 \subset D_\phi$, for if $a\in (M_\phi)_{sa}$ then by the
 Kadison-Schwarz inequality, $a^2=P_{\phi}(a^2)\geq P_{\phi}(a)^2
 =a^2,$ so $P_{\phi}(a^2)= P_{\phi}(a)^2.$ Since $P_\phi
 (M)=M_\phi$ we thus have $P_\phi (M)=P_\phi (D_\phi).$ We have
 therefore shown the following corollary to Theorem ~\ref{thm1}.

 \begin{cor}\label{cor1}
 Let $M$ be a von Neumann algebra and $\phi$ a normal unital map
 in $B(M,M)^+$ such that there exists a faithful normal state
 $\rho$ such that $\rho\circ\phi =\rho.$ Let $P_\phi$ be the
 averaging projection for $\phi.$ Then the following conditions are
 equivalent:\enu{i} $P_\phi$ is entanglement breaking.
 \enu{ii} $M_\phi = P_{\phi}(M)$ is an abelian von Neumann algebra.
 \enu{iii} $\tilde{P_{\phi}}$ is a separable positive functional.
 \end{cor}

  A JC-algebra is called \textit{reversible} if it is
 closed under symmetric products $a_1 a_2... a_n + a_n a_{n-1}...a_1$
 for all $n\in\0N, a_{i}\in A.$ Examples of nonreversible Jordan
 algebras are spin factors, which are the norm closed linear span
 of spin systems, i.e. 1 and self-adjoint unitaries $s_i ,i\in
 I\subset\0N$, satisfying $s_i s_j + s_j s_i =0$ for $i\neq j,$
 where $Card  I\geq 7.$ By \cite{st2}Corollary 7.3, if $P$ is a
 faithful positive projection of a $C^*-$algebra $B$ into itself,
 then $P$ is decomposable if and only if the image  $P(B_{sa})$ of
 $B_{sa}$ is a reversible JC-algebra. We get  the following
 corollary to Corollary ~\ref{cor1}.

 \begin{cor}\label{cor2}
 Let $M$ be a von Neumann algebra and $\phi\in B(M,M)^+$ a unital
 normal map.  Suppose there exists a faithful normal state $\rho$
 such that $\rho\circ\phi =\rho.$  Then if the JW-algebra
 $(M_\phi)_{sa}$ is nonreversible, then $\phi$ is nondecomposable.
 In particular there exists $x\in (A\widehat{\otimes}\1T)^+$ such
 that $(\iota\otimes t)(x)\geq 0,$ while $\tilde{\phi}(x)< 0.$
 \end{cor}

 \bp
 From the above discussion $(M_\phi)_{sa}$ is nonreversible if and
 only if $P_\phi$ is nondecomposable.  If $\phi$ is decomposable,
 so is $\phi ^n$ for all $n\in \0N$, hence $P_\phi$ is
 decomposable, so that $(M_\phi)_{sa}=P_\phi (M_{sa})$ is
 reversible, contrary to assumption.  Thus $\phi$ is
 nondecomposable. By Proposition~\ref{prop1} the existence of $x$
 is clear.

\medskip
 When $M_\phi$ is finite dimensional the structure of $P_\phi$ and
 its dual $\tilde{P_\phi}$ is described in the following more
 general result. Recall that the \textit{centralizer} of a state
 $\rho$ on a $C^*-$algebra $A$ is the set of $a\in A$ such that
 $\rho(ab)=\rho(ba)$for all $b\in A.$

 \begin{prop}\label{prop2}
 Let $A\subset B$ be unital $C^*-$algebras acting on a Hilbert
 space $H$ with $A$  finite dimensional abelian generated by its
 minimal projections $e_1,...,e_n.$ Let $P\in B(B,A)^+$ be a
 unital projection of $B$ onto $A$, and suppose $\rho$ is a
 faithful state on $B$ such that $\rho\circ P =\rho.$ Then we
 have:
 \enu{i} $e_i$ belongs to the centralizer of $\rho$ for all i.
 \enu{ii} If $\rho_{i} (a)=\rho(e_i)^{-1} \rho(e_{i} a e_i), a\in B$,
 then
 $$
 P(a)=\sum_{i=1}^n \rho_{i}(x)e_i.
 $$
 \enu{iii} Let $E\colon B\otimes B(H)\to B\otimes B(H)$ be defined by
 $$
 E(x)=\sum_{i=1}^{n} \rho(e_i)^{-1}(e_{i}\otimes e_{i}^t)x(e_{i}\otimes
 e_{i}^t).
 $$
 Then
 $$
 \tilde{P}=(\rho\otimes Tr)\circ E.
 $$
 \end{prop}

 \bp
 (i) If $a\in B$ then $P(ae_i)=P(a)e_i =e_i P(a)= P(e_i a),$ hence
$$
\rho(e_i a) =\rho(P(e_i a))=\rho(P(ae_i))=\rho(ae_i),
$$
proving (i).

(ii) If $a\in A$ then by (i)
$$
\rho(a) =\sum_{i=1}^n \rho(e_i a)=\sum_{i=1}^n \rho(e_i
ae_i)=\sum_{i=1}^n \rho(P(e_i ae_i))=\sum_{i=1}^n \rho(e_i
P(a)e_i).
$$
Since $A$ is abelian generated by the $e_i's,$ $P(a)=\sum_{i=1}^n
\omega_i(a)e_i$ with $\omega_i$ a state on $B.$ Thus
$$
\rho(a)=\sum_{i=1}^n \rho(e_i(\sum_{j=1}^n\omega_j(a)e_j)e_i)
=\sum_{i=1}^n\omega_i(a)\rho(e_i).
$$
Note that $e_i =P(e_i)=\sum_{j=1}^n\omega_j (e_i)e_j$, hence
$\omega_i (e_i)=1,$ so $\omega_i(e_j)=\delta_{ij}.$ Thus if
 $a=e_j ae_j$ for $i\neq j,$ then $\omega_i(a)=0,$ and
 $\rho(a)=\omega_j(a)\rho(e_j).$
Let $\rho_j$ be as in the statement of the proposition. Then
$\rho_j$ is a state on $B$ with support $e_j,$ and we have
$$
P(a)=\sum_{j=1}^n\omega_j (a)e_j=\sum_{j=1}^n\omega_j (e_j
ae_j)e_j= \sum_{j=1}^n\rho_j(a)e_j,
$$
proving (ii).

 (iii)If $a\otimes b \in B\widehat{\otimes}\1T$ then
\begin{eqnarray*}
\tilde{P}(a\otimes b) &=& Tr(\sum_{i=1}^n\rho_i(a)e_i b^t)\\
&=&\sum_{i=1}^n\rho_i(a)Tr(e_i b^t)\\
&=&\sum_{i=1}^n\rho_i(a)Tr(e_{i}^{t}b)\\
&=&\sum_{i=1}^n\rho(e_i)^{-1} \rho\otimes Tr((e_i\otimes
e_{i}^t)(a\otimes b)(e_i\otimes e_{i}^t))\\
&=& \rho\otimes Tr \circ E.
\end{eqnarray*}
The proof is complete.

\section {Separability in $C^*-$algebras}

The Horodecki Theorem \cite{3 Hor} states that if $M$ and $N$ are
full matrix algebras, and $\rho$ is a state on $M\otimes N$ with
density operator {h}, the $\rho$ is separable if and only if
$\iota\otimes\psi (h)\geq 0$ for all $\psi\in B(N,M)^+.$  We shall
in the present section show a version of this when $M$ is a
nuclear $C^*-$algebra and $N$ a UHF-algebra. In that case we
cannot identify a state with a density operator, so we must
reformulate the theorem. First recall that a $C^*-$algebra $A$ is
\textit{nuclear} if there exists a net of triples
$(M_{n_\lambda}(\0C), \alpha_\lambda, \beta_\lambda),$ where
$\alpha_\lambda\colon A\to M_{n_\lambda}(\0C)$, and $\beta_\lambda
\colon M_{n_\lambda}(\0C)\to A$ are completely positive maps such
that $\lim_\lambda \beta_\lambda\circ \alpha_\lambda (a)=a$ in
norm for all $a\in A,$ see \cite{T}Ch.XV,§1. A $C^*-$algebra $B$
is a UHF-algebra if there is a strictly increasing sequence
$(N_n)_{n\in \0N}$ of $C^*-$subalgebras isomorphic to full matrix
algebras, whose union is norm dense in $B.$ With $A$ and $B$ as
above $A\otimes B$ has a unique $C^*-$norm.

 We say a state $\rho$ on the tensor product $A\otimes B$ of two
$C^*-$algebras is {\textit{weakly separable}} if $\rho$ is a
$w^*-$limit of finite convex sums of product states.

\begin{thm}\label{thm2}
Let $A$ be a nuclear $C^*-$algebra and $B$ a UHF-algebra. Let
$\rho$ be a state on $A\otimes B.$ Then $\rho$ is weakly separable
if and only if $\rho\circ(\iota\otimes\psi)\geq 0 $ for all
$\psi\in B(B,B)^+$.
\end{thm}

We leave it as an open problem how to extend it to more general
$C^*-$algebras. The proof will be divided into some lemmas. One of
them is the finite dimensional version of the theorem. The first
lemma is essentially a restatement of Theorem 3.11 in \cite{st5}.

\begin{lem}\label{lem3}
Let $K$ and $H$ be finite dimensional Hilbert spaces with $ dim
K\leq dim H.$ Let $A\subset B(K)$ be a $C^*-$algebra. Then
$$
A^+ \otimes B(H)^+ =\{x\in A\otimes B(H):
(\iota\otimes\psi)(x)\geq 0, \forall \psi\in B(B(H),H)^+\}
$$
\end{lem}
\bp Let $\1P$ denote the cone on the right side above. Let
$\phi\in B(A,H)^+$ and let $E_A$ be the trace invariant
conditional expectation of $B(K)$ onto $A.$ Then $E_A\circ\phi$ is
a positive extension of $\phi$ to a map in $B(B(K),H)^+.$ Since
$dim K\leq dim H$, $\phi$ is $B^{2}(H)^+$-positive by \cite{st3}
Theorem 3.11, hence $\tilde{\phi}$ is positive on $\1P.$ By
\cite{st3} Lemma 2.1, every positive linear functional on $A\otimes
B(H)$ is of the form $\tilde{\phi}$ with $\phi\in B(A,H)^+$. Thus
every positive linear functional on $A\otimes B(H)$ which is
positive on $A^+ \otimes B(H)^+$ is positive on $\1P$. Since
clearly $A^+\otimes B(H)^+ \subset\1P$ the two cones are equal by
the Hahn-Banach Theorem. The proof is complete.

\begin{lem}\label{lem4}
Let $A,$ $K$ and $H$ be as in the previous lemma, and let $\rho$
be a state on $A\otimes B(H).$ Then $\rho$ is separable if and
only if $\rho\circ (\iota\otimes\psi)\geq 0$ for all $\psi\in
B(B(H),H)^+.$
\end{lem}
\bp If $\rho$ is separable it is clear that $\rho\circ
(\iota\otimes\psi)\geq 0$ for all $\psi\in B(B(H),H)^+.$ To prove
the converse let $Tr$ be the trace on $A\otimes B(H)$ which is 1
on each minimal projection. If $\eta\in B(A\otimes B(H),A\otimes
B(H))$ then its adjoint map $\eta^*$ is defined by
$$
Tr(\eta(x)y)=Tr(x\eta^*(y)).
$$
Then $\eta$ is positive if and only if $\eta^*$ is positive.

Let $h$ be the density matrix for $\rho$ in $ A\otimes B(H)$, and
let $a\in A, b\in B(H).$ Then if $\psi\in B(B(H),H)^+$ we have
$$
\rho\circ(\iota\otimes\psi)(a\otimes
b)=Tr(h(\iota\otimes\psi)(a\otimes
b))=Tr((\iota\otimes\psi^*)(h)(a\otimes b)),
$$
hence
$$
\rho\circ(\iota\otimes\psi)(x)= Tr((\iota\otimes\psi^*)(h)x)
$$
for $x\in A\otimes B(H).$ Thus $(\iota\otimes\psi^*)(h)\geq 0$ for
all $\psi^*\in  B(B(H),H)^+$ and hence for all $\psi\in
B(B(H),H)^+.$ By Lemma 3 $h\in A^+ \otimes B(H)^+,$ so $\rho$ is
separable, completing the proof.

\medskip
The above lemma is false if $dim K > dim H,$ even when $dim K=4$
and $dim H=2$ Indeed, by a result of Woronowicz and the author,
see \cite{W} ,each map in $B(M_2(\0C),M_2(\0C))^+$ is
decomposable. It follows that if $\rho\circ(\iota\otimes t)\geq
0,$ then $\rho\circ(\iota\otimes \psi)\geq 0$ for all $\psi\in
B(M_2(\0C),M_2(\0C))^+,$ hence $\rho$ satisfies the condition of
the lemma. However, there exist maps when $A=M_4 (\0C)$ which
satisfy the Peres (or the PPT) condition which are not separable,
\cite{1Hor}, hence the conclusion of the lemma is false in this
case.

 An easy consequence of the above lemma is the following version
 of the Horodecki Theorem. Note that the difference from the
 Horodecki Theorem is that the maps $\psi$ now map $B(H)$ into
 itself and not into $B(K).$ However, since $dim K\leq dim H$, it
 is easy to deduce the corollary from the Horodecki Theorem.

 \begin{cor}\label{lem3}
 Let $A,$ $K$ and $H$ be as in the above lemma, and let $\rho$ be a state on $A\otimes
 B(H)$with density matrix $h.$ Then $\rho$ is separable if and
 only if $(\iota\otimes \psi)(h)\geq 0$ for all $\psi\in
 B(B(H),H)^+$.
 \end{cor}
 \bp
 Let $Tr$ be as in the above proof. Then
 \begin{eqnarray*}
(\iota\otimes \psi)(h)\geq 0  \ \forall \  \psi \in B(B(H),H)^+
&\Leftrightarrow & Tr((\iota\otimes \psi)(h)x)\geq 0 \  \forall \ x\geq
0\\
&\Leftrightarrow& \rho((\iota\otimes \psi^*)(x)))\geq 0\ \forall\ x,
\psi.
\end{eqnarray*}
By the last lemma this holds if and only if $\rho$ is separable.
The proof is complete.
\medskip

 \textit{Proof of Theorem}.
 If $\rho$ is a state on $A\otimes B$ which is a convex sum of
 product states, then clearly $\rho\circ(\iota\otimes\psi)\geq 0$
 for all $\psi\in B(B,B)^+$. Hence if $\rho$ is a w*-limit of
 convex sums of product states the same inequality holds.

 Conversely assume $\rho\circ(\iota\otimes\psi)\geq 0$ for all $\psi\in
 B(B,B)^+$.  We want to show $\rho$ is a w*-limit of states which are
 convex sums of product states.  Let $x_1,x_2,...,x_n \in A\otimes
 B$ and $\varepsilon > 0$. We must find a state of the form
$\sum_{i=1}^{n} \lambda_i \rho_i\otimes \omega_i$ on $A\otimes B$
such that
$$
|\rho(x_j) - \sum_{i=1}^{n} \lambda_i \rho_{i} \otimes
\omega_i(x_j)|< \varepsilon
$$
for all j. To accomplish this we may assume the $x_j$ belong to
the algebraic tensor product of $A$ and $B$, hence are of the form
$$
x_j=\sum_{i=1}^{n_j} a_{ji}\otimes b_{ji},  \  a_{ji}\in A,
b_{ji}\in B.
$$
Since $B$ is a UHF-algebra there exists a strictly increasing
sequence $N_i$ of subalgebras of $B$ with $N_i\simeq M_{n_i}(\0C)$
such that $B$ is the norm closure of $\bigcup_{i}N_i$.  By
continuity of $\rho$ we may assume there is $k_0$ such  that
$b_{ji}\in N_k$ for $k\geq k_0.$ Let $E_k \colon B\to N_k$ be the
trace invariant conditional expectation of $B$ onto $N_k.$ Thus
each $\psi\in B(N_k ,N_k)^+$ has an extension $\psi\circ E_k$ in
$B(B,B)^+$. Hence $\rho\circ(\iota\otimes\psi)\geq 0$ for all
$\psi\in B(N_k ,N_k)^+$.

Since $A$ is nuclear there exists a net of triples
$(M_{n_{\lambda}}(\0C),\alpha_{\lambda}, \beta{_\lambda})$ such
that $\alpha_{\lambda}\colon A\to M_{n_{\lambda}}(\0C),
\beta_{\lambda}\colon\ M_{n_{\lambda}}(\0C)\to A$ are completely
positive and
$$
\lim_{\lambda} \|\beta_{\lambda}\circ\alpha_{\lambda} (x) -x\| =0
$$
for all $x\in A.$ If we can show $\rho\circ
(\beta_{\lambda}\circ\alpha_{\lambda} \otimes E_k)$ is separable
for arbitrary large $\lambda $ and $k$ such that $n_{\lambda} \leq
n_k$, then $\rho(x_j)$ is arbitrarily well approximated by
separable states for each j. But
$$
\rho\circ(\beta_{\lambda}\circ\alpha_{\lambda} \otimes E_k) =
\rho\circ(\beta_{\lambda} \otimes E_k)\circ (\alpha_{\lambda}
\otimes\iota).
$$
Thus it suffices to show that $\rho\circ(\beta_{\lambda} \otimes
E_k)$ is separable.

We have thus reduced the proof to showing that
$\rho\circ(\beta_\lambda \otimes \iota)$ is a separable functional
on $M_\lambda \otimes N_k$ where $n_k\geq n_\lambda.$ For this let
$y\in (M_\lambda \otimes N_k)^+$. Since $\beta_\lambda$ is
completely positive $\beta_\lambda \otimes \iota(y)\geq 0$ in
$A\otimes N_k$, hence by assumption, if $\psi\in B(N_k,N_k)^+$
$$
\rho((\beta_\lambda \otimes\iota)\circ(\iota\otimes\psi)(y))=
\rho\circ(\iota\otimes\psi)(\beta_\lambda \otimes \iota(y))\geq 0.
$$
It follows that the assumptions of Lemma 10 are satisfied, hence
that $\rho\circ(\beta_\lambda \otimes \iota)$ is separable,
completing the proof of Theorem 8.

\medskip

It is obvious that a convex sum of product states on a tensor
product of $C^*-$algebras can be extended to a similar state on a
tensor product of larger algebras. This is not obvious for weakly
separable states. However, we have

\begin{cor}\label{cor3}
Let $A$ be a unital nuclear $C^*-$algebra and $A_0 \subset A$  a
$C^*-$ subalgebra containing the identity. Let $B$ be a
UHF-algebra. The every weakly separable state on $A_0\otimes B$
has a weakly separable extension to $A\otimes B.$
\end{cor}
\bp Since $A$ is nuclear, so is $A_0.$ Let $\rho$ be a weakly
separable state on $A_0\otimes B$. Let $C$ denote the norm closed
cone in $A\otimes B$ generated by all operators  of the form
$\iota\otimes\psi(x)$ with $x\in (A\otimes B)^+$ and $\psi\in
B(B,B)^+$. Then $C\bigcap (A_0\otimes B)$ is the similar cone with
$x$ now in $A_0 \otimes B$. By the above theorem $\rho$ is
positive on $C\bigcap (A_0\otimes B)$, and by the theorem it
suffices to show that $\rho$ has an extension to $A\otimes B$
which is positive on $C$. Since $\iota \otimes\iota(x)=x$ for all
$x$, $C\bigcap (A_0\otimes B)\supset (A_0\otimes B)^+$, which has
$1\otimes 1$ as an interior point. It follows by a theorem of
Krein-Rutman \cite{Sch}Ch.V,5.4, that $\rho$ has an
 extension which is positive on $C,$
 completing the proof.

Department of Mathematics, University of Oslo, 0316 Oslo, Norway.

e-mail erlings@math.uio.no

\end{document}